\documentclass[a4paper, 12pt]{article}

\usepackage{hyperref}

\usepackage[latin1]{inputenc}
\usepackage[T1]{fontenc}

\usepackage[dvips]{graphicx}
\usepackage{color}

\usepackage{amsthm}  
\usepackage{amsmath} 
\usepackage{amssymb} 
\usepackage{latexsym} 

\theoremstyle{definition}

\newtheorem{proposition}{Proposition}
\newtheorem{observation}{Observation}
\newtheorem{theorem}{Theorem}

\newtheorem{lemma}{Lemma}

\usepackage[numbers]{natbib}
\usepackage{subfigure}

\title{A Remark on Triangle-Critical Graphs}

\author{ \small Anders Sune Pedersen (\texttt{asp@imada.sdu.dk}) \\
  \small Dept. of Mathematics \& Computer Science, University of Southern Denmark, \\
  \small Campusvej 55, 5230 Odense M, Denmark \vspace{5mm} \\
  \small MR Subject Classification: 05C69, 05C05 }

\begin{document}

\maketitle

\begin{abstract}
  A connected $k$-chromatic graph $G$ with $k \geq 3$ is said to be
  triangle-critical, if every edge of $G$ is contained in an induced
  triangle of $G$ and the removal of any triangle from $G$ decreases
  the chromatic number of $G$ by three. B.  Toft posed the problem of
  showing that the complete graphs on more than two vertices are the
  only triangle-critical graphs. By applying a method of M. Stiebitz
  [Discrete Math. 64 (1987), 91--93], we answer the problem
  affirmatively for triangle-critical $k$-chromatic graphs with $k
  \leq 6$.
\end{abstract}

\section{Introduction}
The purpose of this manuscript is to introduce the Triangle-Critical
Graph Conjecture, which is a special case of Lovasz' Double-Critical
Graph Conjecture, and settle that conjecture for all $k$-chromatic
graphs with $k \leq 6$. For $k \leq 5$, the desired result follows
directly from a theorem of M. Stiebitz \cite{MR882614}. We show that
the method of M.  Stiebitz \cite{MR882614} extends so as to settle the
Triangle-Critical Graph Conjecture for $k \leq 6$. All graphs
considered in this manuscript are assumed to be simple and finite. The
reader is referred to \cite{DiestelGT2006} for definitions of any
graph-theoretic concept used but not explicitly defined in this
manuscript. A graph $G$ is called \emph{vertex-critical} if $\chi
(G-v) < \chi(G)$ for every vertex $v \in V(G)$. Clearly, every
vertex-critical graph is connected. A critical graph $G$ is called
\emph{double-critical} if the chromatic number of $G$ decreases by at
least two whenever two adjacent vertices are removed from $G$, that
is, for any edge $xy \in E(G)$, we have $\chi (G - x - y ) = \chi(G) -
2$. Clearly, all complete graphs are double-critical, and L. Lovasz
\cite{MR0232693} conjectured that only the complete graphs are
double-critical. M. Stiebitz \cite{MR882614} proved the
Double-Critical Graph Conjecture for all double-critical $k$-chromatic
graphs with $k \leq 5$, while the conjecture remains open for $k \geq
6$. Given the difficulties in settling the Double-Critical Graph
Conjecture it might be enlightening to study weaker conjectures on the
same theme. For instance, if we require the chromatic number $\chi
(G)$ to decrease by two whenever we remove two independent edges from
$G$, then it is straightforward to see that $G$ is a complete graph.
Similarly, if we consider removing an arbitrary edge $uv$ and an
arbitrary vertex $w \notin \{ u, v \}$, and require the chromatic
number to decrease by two. The aforementioned weakenings of the
Double-Critical Graph Conjecture are of little interest, since they
are straightforward to prove and the proof ideas do not extend to
double-critical graphs in generel. In the following we study a
weakening of the Double-Critical Graph Conjecture, which does not seem
to be readily solvable.

At the graph theory meeting GT2007 held at Fredericia, Denmark, on the
sixth to the ninth of December 2007, B. Toft posed the following
weakening of Lovasz' Double-Critical Graph Conjecture: \newline
\newline Given a $k$-chromatic ($k \geq 3$) connected graph $G$ with
the property that any edge $e$ is contained in a triangle, and for any
three vertices $x,y$ and $z$ which induce a triangle, $\chi (G - x - y
- z) = k - 3$ (we shall refer to such graphs as
\emph{triangle-critical}). Is it true that $G$ is the complete
$k$-graph?  \newline \newline

We shall refer to this special case of the Double-Critical Graph
Conjecture as the \emph{Triangle-Critical Graph Conjecture}. Clearly,
any triangle-critical graph is vertex-critical and double-critical,
and so, by theorem of M. Stiebitz \cite{MR882614}, the
Triangle-Critical Graph Conjecture is settled affirmatively for every
$k \leq 5$. In Section \ref{sec:sixChromatic}, we apply the method of
M. Stiebitz \cite{MR882614} to show that the only $6$-chromatic
triangle-critical graph is the complete $6$-graph.  However, the
Triangle-Critical Graph Conjecture does not seem to be substantially
easier to settle for graphs with chromatic number greater than six. In
Section \ref{sec:sixChromatic}, we shall need the following result.
\begin{proposition}
  Suppose $G$ is a double-critical $k$-chromatic graph, which is not a
  complete graph. Then $G$ does not contain a complete $(k-1)$-graph
  as a subgraph.
\label{prop:forbiddenCompleteKminusOne}
\end{proposition}
The proof of Proposition~\ref{prop:forbiddenCompleteKminusOne} is
elementary and, therefore, omitted.

\section{On triangle-critical $6$-chromatic graphs}
\label{sec:sixChromatic}
Let $[ k ]$ denote the set $\{1,2, \ldots, k \}$. Any $k$-colouring
will, unless otherwise stated, by assumed to employ the colours of the
set $[k]$. Given a colouring $c$ of some graph $G$ and a subset $S$ of
the vertex-set $V(G)$, let $c(S)$ denote the set of colours that $c$
apply to the vertices of $S$. Two colourings $c_1$ and $c_2$ of $G$
are \emph{equivalent} if, for all vertices $u,v \in V(G)$, $c_1(u) =
c_1(v)$ if and only if $c_2(u) = c_2(v)$. A graph $G$ with $\chi (G)
\leq 4$ is said to be \emph{uniquely $4$-colourable} if any two proper
$4$-colourings of $H$ are equivalent. Recall, that any $k$-colouring
may also be considered as a $j$-colouring for any integer $j > k$.

\begin{observation}
\begin{itemize}
\item[]
\item[(a)]
If $H$ is a uniquely $4$-colourable graph with $\chi (H) \leq 3$, then
$H \in \{ K_1, K_2, K_3 \}$.
\item[(b)]
If $H$ is a uniquely $4$-colourable graph with $n(H) = 4$, then
$H \simeq K_4$.
\end{itemize}
\label{obs:uniqueFourColourings}
\end{observation}
\begin{proof}
\begin{itemize}
\item[(a)] Suppose $H$ is a uniquely $4$-colourable graph with $\chi
  (H) \leq 3$. Moreover, suppose that there exists a pair of
  non-adjacent vertices in $H$, say $u,v \in H$ with $uv \notin E(H)$.
  Let $c$ denote some $3$-colouring of $H$, and define another proper
  colouring $c'$ of $H$ as follows. For every vertex $w \in V(H)
  \backslash \{ u, v \}$, define $c'(w) := c(w)$. If $c(u) = c(v)$,
  define $c'(v) := c(v)$ and $c'(u) := 4$; otherwise, if $c(u) \neq
  c(v)$, define $c'(u) = c'(v) := 4$. In any case, $c$ and $c'$ are
  two non-equivalent proper $4$-colouring of $H$, which contradiction
  the assumption that $H$ is uniquely $4$-colourable. Thus, $H$ must
  be complete, and, since $\chi (H) \leq 3$, we obtain $H \in \{ K_1, K_2, K_3
  \}$.
\item[(b)] Let $H$ denote a uniquely $4$-colourable graph with $n(H) =
  4$. If $\chi (H) \leq 3$, then (a) implies $n(H) \leq 3$. Thus, we
  must have $\chi (H) \geq 4$, and, since $H$ is $4$-colourable,
  indeed, $\chi (H) = 4$. Obviously, the only $4$-chromatic on four
  vertices is $K_4$.
\end{itemize}
\end{proof}


We denote the number of colours that $c$ applies to the vertices of
$S$ by $|c(S)|$. Given some subset $S$ of $V(G)$, let $T(S:G)$ denote
the common neighbours of the vertices of $S$ in $G$.

In the following we let $G$ denote an arbitrary $6$-chromatic
triangle-critical graph.

\begin{lemma}
  Suppose that $x, y$ and $z$ are three distinct vertices of $G$ such
  that $G[x,y,z] \simeq K_3$. Then $|c(T(x,y,z : G))| = 3$ for any
  proper $3$-colouring of $G - \{x,y,z\}$. In particular, $|T(x,y,z :
  G)| \geq 3$.
\label{lem:SizeOfCommonNeighbourhood}
\end{lemma}

\begin{proof}
  Since $G$ is $6$-chromatic triangle-critical and $G[x,y,z] \simeq
  K_3$, we obtain $\chi (G - x - y - z ) = 3$. Let $c$ denote an
  arbitrary proper $3$-colouring of $G - x - y - z$. Obviously, the
  colouring $c$ can be extended to a proper $6$-colouring of $G$
  defining $c(x)=4$, $c(y)=5$ and $c(z)=6$. It is easy to see that in
  any colour-class $S_i$ for $i \in [6]$ there must be a vertex $z_i$
  adjacent to at least one vertex in each of the other five
  colour-classees of $G$ under the colouring $c$; otherwise, the
  vertices of some colour-class $S_i$ could be recoloured using the
  colours $[6] \backslash \{ i \}$ to produce a proper $5$-colouring
  of $G$. Now the three vertices $z_1, z_2$ and $z_3$ are all adjacent
  to the vertices $x,y$ and $z$, that is, $\{ z_1, z_2, z_3 \}
  \subseteq T(x,y,z : G)$ and $|c(T(x,y,z : G))|=3$.
\end{proof}

\begin{theorem}
The only $6$-chromatic triangle-critical graph is the the complete $6$-graph.
\label{th:SixChrom}
\end{theorem}

The proof of Theorem \ref{th:SixChrom} follows the lines of an
argument by M. Stiebitz \cite{MR882614}.

\begin{proof}
  Let $G$ denote an arbitrary $6$-chromatic triangle-critical graph,
  and suppose that $G$ is not the complete $6$-graph. Let $\{ v_1,
  v_2, \ldots, v_n \}$ denote the vertex set of $G$, and let $H_1,
  H_2, \ldots, H_r$ denote a sequence of induced subgraphs of $G$ such
  that
\begin{itemize}
\item[(i)]
for every $i \in \{ 1, \ldots, r \}$, the graph $H_i$ is uniquely
$4$-colourable;
\item[(ii)]
for every $i \in \{ 1, \ldots, r \}$, the graph $H_i$ is of order $i$
and $V(H_i) = \{ v_1, \ldots, v_i \}$;
\item[(iii)]
for every $i \in \{ 1, \ldots, r-1 \}$, the graph $H_i$ is a subgraph of
$H_{i+1}$; and
\item[(iv)]
there is no uniquely $4$-colourable induced subgraph of $G$ of order
$r+1$ containing $H_r$ as a subgraph.
\end{itemize}
Now, the graph $G$ contains an edge and so, by definition, $G$ contains at least one triangle, say $G[x,y,x] \simeq K_3$. By
Lemma~\ref{lem:SizeOfCommonNeighbourhood}, there exists some vertex $w
\in T(x,y,z:G)$. Now $G[x,y,x,w] \simeq K_4$, which is a uniquely
$4$-colourable induced subgraph of $G$ of order $4$, and so $r \geq
4$. Obviously, $H_r \neq G$, since $\chi(H_r) = 4$ and $\chi (G) = 6$.

Since $r \geq 4$, it follows from
Observation~\ref{obs:uniqueFourColourings}~(b), that $H_4 \simeq K_4$ and
so $G[v_2, v_3, v_4] \simeq K_3$. In particular, $H_r$ contains $K_3$
as a subgraph. Let $v_i, v_j$ and $v_l$ denote three vertices of
$H_r$ with say $i > j > l$ such that $G[v_i, v_j, v_l ] \simeq K_3$
and suppose that the vertices $v_i$, $v_j$ and $v_l$ are choosen such
that $l$ is maximal. Since $G[v_4, v_3, v_2] \simeq K_3$, $l \geq
2$. Now we claim that $T(v_i, v_j, v_l : H_r) \subseteq
V(H_{l-1})$. If the common neighbourhood $T(v_i, v_j, v_l : H_r)$ were
not a subset of $V(H_{l-1})$, then there would exist some vertex
$v_q \in T(v_i, v_j, v_l : H_r)$ with $q > l$. But then, since
$G[v_i,v_j,v_q] \simeq K_3$, we would have a contradiction with the
maximality of $l$.

Let $c_3$ denote an arbitrary proper $3$-colouring of $G \backslash \{
v_i, v_j, v_l \}$. Since $H_{l-1}$ is a subgraph of $G \backslash \{
v_i, v_j, v_l \}$, the graph $H_{l-1}$ is $3$-colourable. Now, it
follows from Observation~\ref{obs:uniqueFourColourings}~(a), that
$H_{l-1}$ is a complete graph (on less than four vertices).  Since
$T(v_i, v_j, v_l : H_r) \subseteq V(H_{l-1})$, it follows that the
vertices of $T(v_i, v_j, v_l : H_r)$ induce a complete subgraph of
$G$.

Since $H_r$ is $4$-colourable, it follows that the vertices of $T(v_i,
v_j, v_l : H_r)$ must be non-adjacent in $G$. Thus, the induced graph
$G[ T(v_i, v_j, v_l : H_r) ]$ is both complete and edge-empty, which
implies that the set $T(v_i, v_j, v_l : H_r)$ consists of at most one
vertex.

According to Lemma~\ref{lem:SizeOfCommonNeighbourhood}, $|T(v_i, v_j,
v_l : G)| \geq 3$. Thus, since $|T(v_i, v_j, v_l : H_r)| \leq 1$, the
set $T(v_i, v_j, v_l : G) \backslash T(v_i, v_j, v_l : H_r)$ contains
at least two distinct vertices, say $u$ and $v$. We shall consider the
graphs $H_u := G[V(H_r) \cup \{ u \}]$ and $H_v := G[V(H_r) \cup \{ v
\}]$. The maximality of $r$ implies that either (i) $H_u$ is
non-uniquely $4$-colourable or (ii) $H_u$ is not $4$-colourable, i.e.,
$\chi (H_u) = 5$. It is easily seen that if $H_u$ were $4$-colourable,
then it would be uniquely $4$-colourable. Thus, it must be the case
that $\chi ( H_u ) = 5$ and, similarly, $\chi (H_v) = 5$. Let $S :=
V(G) \backslash V(H_r)$. If there were an edge in $S \backslash \{ u
\}$ or $S \backslash \{ v \}$, then, since $G$ is also
double-critical, we would have $\chi (H_u) \leq 4$ or $\chi (H_v) \leq
4$, a contradiction. On the other hand, if $G[S]$ were edge-empty,
then the proper $4$-colouring of $H_r$ could be extended to a proper
$5$-colouring of $G$, a contradiction. Thus, it must be the case that
$u$ and $v$ are adjacent in $G$, and so it follows that the set $\{
v_i, v_j, v_l, u, v \}$ induce the complete $5$-graph in $G$. Now we
have a contradiction with
Proposition~\ref{prop:forbiddenCompleteKminusOne}, and so it must be
the case that $G$ is the complete $6$-graph.
\end{proof}

\subsection*{Acknowledgement}
The author wish to thank B. Toft for useful discussions on
double-critical and triangle-critical graphs.

\bibliographystyle{plainnat}  
\bibliography{natbibFull}

\end{document}